\documentclass [a4paper,twoside,11pt]{article}

\usepackage[latin1]{inputenc}
\usepackage{amsmath,amsfonts,amssymb}
\usepackage{vmargin,graphicx,theorem}
\usepackage[english]{babel}
\usepackage{enumerate}
\usepackage{setspace}
\usepackage{color}
\setpapersize[portrait]{A4}
\setmarginsrb{2.7cm}{2cm}{2.7cm}{2cm}{0cm}{0.1cm}{0.5cm}{2cm}
\newcommand{\disp}{\displaystyle}

\newcommand{\dN}{\ensuremath{\mathbb{N}}}

\newcommand{\dR}{\ensuremath{\mathbb{R}}}

\newcommand{\dZ}{\ensuremath{\mathbb{Z}}}

\newtheorem{ethm}{Theorem}[section]

\newtheorem{ecor}[ethm]{Corollary}

\newtheorem{eprop}[ethm]{Proposition}

\newtheorem{elem}[ethm]{Lemma}

\newtheorem{edefi}[ethm]{Definition}

\newtheorem{exe}[ethm]{Example}

\newcommand{\proofend}{~$\rhd$}
\newcommand{\proofbegin}{~$\lhd$}

\newenvironment{eproof}
 {\noindent {\emph{\textbf{Proof}}}\proofbegin~}
 {\proofend}

\newcommand{\p}[4]{{#3}\!\left#1{#4}\right#2}

\newcommand{\ABS}[1]{\ensuremath{{\left| #1 \right|}}} 
\newcommand{\PAR}[1]{\ensuremath{{\left(#1\right)}}} 
\newcommand{\SBRA}[1]{\ensuremath{{\left[#1\right]}}} 
\newcommand{\BRA}[1]{\ensuremath{{\left\{#1\right\}}}} 
\newcommand{\NRM}[1]{\ensuremath{{\left\Vert #1\right\Vert}}} 

\renewcommand{\phi}{\varphi}

\renewcommand{\leq}{\leqslant}
\renewcommand{\geq}{\geqslant}

\newcommand{\varf}[1]{\mathbf{Var}_{#1}}
\newcommand{\entf}[1]{\mathbf{Ent}_{#1}}

\newcommand{\ent}[2]{\p(){\entf{#1}}{#2}}

\newcommand{\var}[2]{\p(){\varf{#1}}{#2}}

\def\disp{\displaystyle}

\newcommand{\R}{\dR}

\newcommand{\C}[1]{\ensuremath{{\mathcal C}^{#1}}}
\newcommand{\Cn}[1]{\mathcal{C}^{#1}}

\newcommand{\GI}{{\mathcal L}}

\newcommand{\al}{\alpha}

\newcommand{\be}{\beta}

\newcommand{\ep}{\epsilon}

\newcommand{\1}{\mathbb I}

\newcommand{\osc}[1]{{\bf Osc}_\mu\PAR{#1}} 

\newcommand{\email}[1]{{E-mail: {\textsf {#1}}}}

\newcommand{\Capa}{{\rm Cap}}
\newcommand{\WPME}{{\rm (WPME)} }
\newcommand{\CP}{{\rm C}_{\rm P} }
\newcommand{\CLSI}{{\rm C}_{\rm LS} }
\newcommand{\CWL}{{\rm h}_{\rm WLS} }
\newcommand{\CWP}{{\rm \beta}_{\rm WP} }
\newcommand{\COP}{{\rm C}_{\rm OP} }
\newcommand{\betaP}{\beta_{\kern 1pt\rm P} }
\newcommand{\betaSL}{\beta_{\kern 1pt\rm LS} }
\newcommand{\kappaP}{\kappa_{\rm P} }
\newcommand{\kappaLS}{\kappa_{\rm LS} }
\newcommand{\kkappa}{\mbox{\sc k}}
\renewcommand{\(}{\left(}
\renewcommand{\)}{\right)}

\begin{document}
\title{\sl $L^q$-functional inequalities and weighted porous media equations}
\author{Jean Dolbeault, Ivan Gentil, Arnaud Guillin and Feng-Yu Wang}
\date{\today}
\maketitle
\thispagestyle{empty}

\begin{abstract}\noindent
Using measure-capacity inequalities we study new functional inequalities, namely $L^q$-Poincar\'e inequalities
$$
\var{\mu}{f^q}^{1/q}\le \CP\int |\nabla f|^2\,d\mu
$$
and $L^q$-logarithmic Sobolev inequalities
$$
\ent{\mu}{f^{2q}}^{1/q}\le \CLSI \int |\nabla f|^2\,d\mu\;.
$$
As a consequence, we establish the asymptotic behavior of the solutions to the so-called weighted porous media equation
$$
\frac{\partial u}{\partial t} =\Delta u^m-\nabla\psi\cdot\nabla u^m
$$
for $m\ge1$, in terms of $L^2$-norms and entropies.
\end{abstract}


\section{Introduction}
\label{Intro}

In this paper we analyze decay rates of the entropies associated to nonlinear diffusion equations using inequalities relating entropy and entropy production functionals. Consider for instance the Ornstein-Uhlenbeck semi-group on $\dR^d$, which is governed by
$$
\label{eq-ou}
\frac{\partial u}{\partial t} =\Delta u-x\cdot\nabla u\;,\quad t\geq0\;,\quad x\in \R^d\;,
$$
with an initial condition $u_0\in\C{2}\PAR{\R^d}\cap L^1_+(\R^d,d\gamma)$. Here $d\gamma=\gamma\,dx$ is the Gaussian measure on $\dR^d$, $\gamma(x):=(2\pi)^{-d/2}\,\exp(-\ABS{x}^2/2)$. Two entropies are widely used, namely
$$
\var{\gamma}{u}:=\int\PAR{u-\int u\;d\gamma}^2\,d\gamma\quad\text{and}\quad\ent{\gamma}u:=\int u\log\PAR{\frac u{\int u\;d\gamma}}\;d\gamma\;.
$$
If $u$ is a smooth solution of the Ornstein-Uhlenbeck equation, integrations by parts show that
$$
\frac d{dt}\,\var{\gamma}{u}=-\;2\int|\,\nabla u\,|^2\,d\gamma\quad\text{and}\quad\frac d{dt}\,\ent{\gamma}{u}=-\;4\int|\,\nabla \sqrt u\,|^2\,d\gamma\;.
$$
By the Poincar\'e inequality in the first case,
\begin{equation*}
\forall\;f\in\C{1}(\R^d)\;,\quad \var{\gamma}{f}\leq\int\ABS{\nabla f}^2\,d\gamma\;,
\end{equation*}
and Gross' logarithmic Sobolev inequality, see \cite{gross}, in the second case,
$$
\label{eq-logsob}
\forall\;f\in\C{1}(\R^d)\;,\quad \ent{\gamma}{f^2} \leq\; 2\int\ABS{\nabla f}^2\,d\gamma\;,
$$
it follows from Gronwall's lemma that for any $t\geq 0$,
$$
\var{\gamma}{u}\leq e^{-2t}\,\var{\gamma}{u_0}\quad\text{and}\quad\ent{\gamma}{u}\leq e^{-2t}\,\ent{\gamma}{u_0}
$$
for all smooth initial conditions $u_0$ in $L^2_{\gamma}(\dR^d)$ in the first case, or such that $\ent{\gamma}{u_0}$ is finite in the second case. With minor changes, the method can be extended to the semi-group generated by
$$
\label{eq-oug}
\frac{\partial u}{\partial t} =\Delta u-\nabla\psi\cdot\nabla u\;,\quad t\geq0\;,\quad x\in \R^d\;,
$$
if $\psi$ is a smooth function such that $Z:=\int e^{-\psi}\,dx$ is finite and the probability measure $d\mu_\psi=Z_\psi^{-1}\,e^{-\psi}\,dx$ satisfies a Poincar\'e inequality or a logarithmic Sobolev inequality. See for example \cite{logsob} for a review on this inequalities and there applications.

\medskip
A natural question is how to extend the variance or the entropy convergence to nonlinear semi-groups. Let $m>1$, and consider the semi-group generated by the {\it weighted porous media equation\/}
$$
\label{eq-pdein}
\frac{\partial u}{\partial t}=\Delta u^m-\nabla \psi\cdot\nabla u^m\,,\quad t\geq0\;,\quad x\in\dR^d\;,
$$
with a non-negative initial condition $u(0,x)=u_0(x)$. This equation, in short (WPME) is a simple extension of the standard porous media equation, which corresponds to $\psi=0$. We shall refer to \cite{vazquez} for an introduction on this topic. A major difference is that under appropriate conditions on $\psi$ the solution of \WPME converges to its mean. In other words, the nonlinear semi-group converges to the limit measure $\mu_\psi$. The variance of a solution solution of \WPME now obeys to
$$
\frac d{dt}\var{\mu_\psi}{u}=-\,\frac 8{(m+1)^2}\int|\,\nabla u^{\frac{m+1}2}\,|^2\,d\mu_\psi\;.
$$
Classical Poincar\'e and logarithmic Sobolev inequalities are of no more use and have to be replaced by adapted functional inequalities, which are the purpose of this paper. This paper extends some earlier results on solutions to the porous media equation on the torus $S^1\equiv [0,1)$ 	and related functional inequalities, see \cite{car-al}. As for the functional inequalities, we will work in a more general framework involving two Borel probability measures $\mu$ and $\nu$ on a Riemannian manifold $(M,g)$, which are not necessarily absolutely continuous with respect to the volume measure. To consider quantities like $\int f^q\,d\mu$ and $\int\ABS{\nabla f}^2\,d\nu$, it is therefore natural to work in the space of functions $f\in\Cn{1}(M)$, although slightly more general function spaces can be introduced by density with respect to appropriate norms. If the measures were absolutely continuous with respect to the volume measure, we could take functions which are only locally Lipschitz continuous as in \cite{ca-ba-ro2}.

\medskip
In Section~\ref{sec-fi} we will define functional inequalities that we shall call {\it $L^q$-Poincar\'e\/} and {\it $L^q$-logarithmic Sobolev inequalities.\/} Equivalence of these inequalities with capacity-measure criteria will be established, based on Maz'ja's theory. Links with more classical inequalities such as weak Poincar\'e or weak logarithmic Sobolev inequalities are then studied in Section~\ref{Further}. Explicit criteria can be deduced from earlier works, mainly \cite{barthe-roberto,ca-ba-ro2,ca-ge-gu}. In Section~\ref{sec-app} we will give applications to the weighted porous media equation. Using the $L^q$-Poincar\'e and $L^q$-logarithmic Sobolev inequalities we describe the asymptotic behavior of the solutions in terms of variance or entropy. The proof of two variants of results of \cite{ca-ba-ro2} is given in an appendix, see Section~\ref{Sec:Appendix}.

\medskip
Throughout this paper, we intend to work under minimal assumptions and do not require that the measures showing up on both sides of the inequalities are the same or that they are absolutely continuous with respect to the volume measure. However, when only one measure is specified, one has to understand that the measures $\mu$ and $\nu$ are the same on both sides of the inequalities.

\section{Two $L^q$-functional inequalities}
\label{sec-fi}

\subsection{$L^q$-Poincar\'e inequalities}
\label{sec-poincare}

\begin{edefi}
Let $\mu$ and $\nu$ be respectively a probability measure and a positive measure on~$M$. Assume that $q\in(0,1]$. We shall say that $(\mu,\nu)$ satisfies a {\rm $L^q$-Poincar\'e inequality} with constant $\CP$ if for all non-negative functions $f\in\Cn{1}(M)$ one has
\begin{equation}
\label{eq1}
\big[\var{\mu}{f^q}\big]^{1/q}:=\left[{\int f^{2q}\,d\mu-\PAR{\int f^q\,d\mu}^2}\right]^{1/q}\leq\CP\int\ABS{\nabla f}^2\,d\nu\;.
\end{equation}
\end{edefi}

Note that if $q>1$, Inequality~\eqref{eq1} is not true if $\mu$ is not a Dirac measure. Consider indeed $f=1+\ep\,g$ with $\ep\to 0$ and $g$ bounded. By applying Inequality~\eqref{eq1} we get
$$
\Big[{q^2\ep^2\PAR{\var{\mu}{g^2}+o(1)}}\Big]^{1/q}\leq\ep^2\,\CP\int\ABS{\nabla g}^2\,d\nu+o(\ep^2)\;.
$$
If $g$ is such that $\var{\mu}{g^2}$ and $\int\ABS{\nabla g}^2\,d\nu$ are both positive and finite, we obtain a contradiction by letting $\ep\to 0$ if $q>1$.
\begin{eprop}
\label{prop-hi}
For any bounded positive function $f$, the function $q\mapsto\var{\mu}{f^q}^{1/q}$ is increasing with respect to $q\in(0,1]$. As a consequence, if the $L^{q_1}$-Poincar\'e inequality holds, then the $L^{q_2}$-Poincar\'e inequality also holds for any $0<q_2\leq q_1\leq 1$.
\end{eprop}
We shall say that $L^q$-Poincar\'e inequalities form a {\sl hierarchy\/} of inequalities. The classical Poincar\'e inequality corresponding to $q=1$ implies all $L^q$-Poincar\'e inequalities for $q\in(0,1)$. \smallskip

\begin{eproof}
Without loss of generality, we may assume that $f$ is positive. For any $q\in(0,1)$, let $F(q):=\var{\mu}{f^q}^{1/q}$. We have
\begin{equation*}
\label{N1}
\frac{F'(q)}{F(q)}= \frac{d}{d q}\log F(q) =\frac 1q\,\bigg\{-\log F(q) + \frac{\mu(f^{2q}\log f^2)-
\mu(f^q)\,\mu(f^q\log f^2)}{F(q)^q}\bigg\}\;.
\end{equation*}
Let $h(t):=\mu(f^{tq}\log f^2)/\mu(f^{tq})$, $t\in(1,2)$ and observe that by the Cauchy-Schwarz inequality, $\frac 2q\,[\mu(f^{tq})]^2\,h'(t)=\mu(f^{tq}\,(\log f^2)^2)\,\mu(f^{tq})-[\mu(f^{tq}\,\log f^2)]^2\geq 0$. This proves that
$$
h(2)=\frac{\mu(f^{2q}\log f^2)}{\mu(f^{2q})}\geq \frac{\mu(f^q\log f^2)}{\mu(f^q)}=h(1)\;.
$$
Hence we arrive at
\begin{equation*}
\begin{split}
\frac{qF'(q)}{F(q)}\ge -\log F(q) + \frac{\left(1-\frac{[\mu(f^q)]^2}{\mu(f^{2q})}\right)\,\mu(f^{2q}\log f^2)}{F(q)^q}&=-\log F(q) + \frac{\mu(f^{2q}\log f^{2q})}{q\,\mu(f^{2q})}\\
&\quad\ge -\log F(q) +\log \PAR{\mu(f^{2q})}^{1/q}\ge 0
\end{split}
\end{equation*}
where the last two inequalities hold by Jensen's inequality and by monotonicity of the logarithm.\end{eproof}

\medskip
We will now give a characterization of the $L^q$-Poincar\'e inequality in terms of the {\sl capacity measure criterion.\/} Such a criterion has recently been applied in \cite{ca-ba-ro,chen05,ca-ge-gu} to give necessary and sufficient conditions for the usual, weak or super Poincar\'e inequality, and the usual or weak logarithmic Sobolev inequality or the $F$-Sobolev inequality. The capacity measure criterion allows to compare all these inequalities and can be characterized in terms of Hardy's inequality, in the one-dimensional case.

\medskip Let $\mu$ and $\nu$ be respectively a probability measure and a positive measure on~$M$. Given measurable sets $A$ and $\Omega$ such that $A\subset\Omega\subset M$, {\sl the capacity\/} $\Capa_\nu(A,\Omega)$ is defined as
$$
\Capa_\nu(A,\Omega):=\inf\BRA{\int\ABS{\nabla f}^2\,d\nu\;:\;f\in{\mathcal C}^1(M)\,,\;\1_A\leq f\leq\1_\Omega}\;.
$$
If the set $\BRA{f\in{\mathcal C}^1(M)\,:\,\1_A\leq f\leq\1_\Omega}$ is empty then, by convention, we set $\Capa_\nu(A,\Omega):=+\infty$. This the case of $\Capa_\nu(A,A)=+\infty$ for any bounded measurable set $A$ and any $\nu$ with a locally positive density.

\medskip Let $q\in(0,1)$ and define
\begin{equation*}
\label{eq-beta}
\betaP :=\sup\BRA{\;\sum_{k\in\dZ}\frac{\big[\mu(\Omega_k)\big]^{1/(1-q)}} {\big[\Capa_\nu(\Omega_k,\Omega_{k+1})\big]^{q/(1-q)}}\;}^{(1-q)/q}
\end{equation*}
where the supremum is taken over all $\Omega\subset M$ with $\mu(\Omega)\le 1/2$ and all sequences
$\PAR{\Omega_k}_{k\in\dZ}$ such that for all $k\in\dZ$, $\Omega_k\subset\Omega_{k+1}\subset \Omega$.
\begin{ethm}
\label{thm-1}
Let $\mu$ and $\nu$ be respectively a probability measure and a positive measure on~$M$.
\begin{enumerate}
\item [{\rm (i)}] If $q\in[1/2,1)$ and $(\mu,\nu)$ satisfies a $L^q$-Poincar\'e inequality with a constant $\CP$, then $\betaP\leq 2^{1/q}\,\CP$.
\item [{\rm (ii)}] If $q\in(0,1)$ and $\betaP<+\infty$, then $(\mu,\nu)$ satisfies a $L^q$-Poincar\'e inequality with constant $\CP\leq \kappaP\,\betaP $, for some constant $\kappaP$ which depends only on $q$.
\end{enumerate}
\end{ethm}
\begin{eproof}
The proof follows the main lines of Theorem~2.3.5 of \cite{mazja}.

\medskip {\sl Proof of \/}(i). Consider $\Omega\subset M$ such that $\mu(\Omega)\le 1/2$ and let $\PAR{\Omega_k}_{k\in\dZ}$ be a sequence such that for all $k\in\dZ$, $\Omega_k\subset\Omega_{k+1}\subset\Omega$. Fix $N\in\dN^*$ and for $k\in\BRA{-N,\ldots N}$, let $f_k\in{\mathcal C}^1(M)$ be such that $\1_{\Omega_k}\leq f_k\leq\1_{\Omega_{k+1}}$. If no such $f_k$ exists, that is if $\Capa_{\nu}(\Omega_k,\Omega_{k+1})=+\infty$, then we discard $\Omega_{k+1}$ from the sequence and reindex it. Finally, let $(\tau_k)_{k\in\BRA{-N,\ldots N}}$ be a non-increasing family of non-negative reals numbers to be defined later. A function $f$ on $M$ is defined as follows:
\begin{enumerate}
\item[(1)] $f=\tau_{-N}$ on $\Omega_{-N}$,
\item[(2)] $f=(\tau_k-\tau_{k+1})\,f_k+\tau_{k+1}$ on $\Omega_{k+1}\setminus\Omega_k$ for all $k\in\BRA{-N,N-1}$,
\item[(3)] $f=\tau_N\,f_{N+1}$ on $\Omega_{N+1}\setminus\Omega_N$ and $f=0$ on $\Omega\setminus\Omega_{N+1}$.
\end{enumerate}
Using the fact that $f=0$ on $\Omega^c$, it follows from the Cauchy-Schwarz inequality that
$$
\PAR{\int f^q\,d\mu}^2\leq\mu(\Omega)\int f^{2q}\,d\mu\;,
$$
from which we get
$$
\var{\mu}{f^q}\geq \frac{1}{2}\int f^{2q}\,d\mu\;.
$$
By the co-area formula, we obtain
$$
\int f^{2q}\,d\mu=\kern -3pt\int_0^\infty\kern -10pt\mu\PAR{\{f\geq t\}}\,d(t^{2q}) \geq \kern -5pt\sum_{k=-N}^{N-1}\int_{\tau_{k+1}}^{\tau_k}\kern -10pt\mu(\{f\geq\tau_k\})\,d(t^{2q}) = \kern -3pt\sum_{k=-N}^{N-1}\kern -4pt\mu(\{f\geq\tau_k\})\PAR{\tau_k^{2q}-\tau_{k+1}^{2q}}\,.
$$
{}From $2q\geq 1$, we get $\PAR{\tau_k^{2q}-\tau_{k+1}^{2q}}\geq\PAR{\tau_k-\tau_{k+1}}^{2q}$, and
$$
\var{\mu}{f^q}\geq \frac 12\sum_{k=-N}^{N-1}\mu(\Omega_k)\,\PAR{\tau_k-\tau_{k+1}}^{2q}\;.
$$
Using the $L^q$-Poincar\'e inequality we get
$$
\PAR{\,\frac 12\sum_{k=-N}^{N-1}\mu(\Omega_k)\,\PAR{\tau_k-\tau_{k+1}}^{2q}}^{1/q}\leq \CP\int\ABS{\nabla f}^2\,d\nu\;.
$$
On the other hand, with the convention $\tau_{N+1}=0$, we have
$$
\int\ABS{\nabla f}^2\,d\nu=\sum_{k=-N}^N\PAR{\tau_k-\tau_{k+1}}^{2} \int_{\Omega_{k+1}\setminus\Omega_k} \ABS{\nabla f_k}^2\,d\nu\;.
$$
We may now take the infimum over all functions $f_k$ and obtain
$$
\PAR{\,\frac 12\sum_{k=-N}^{N-1}\mu(\Omega_k)\,\PAR{\tau_k-\tau_{k+1}}^{2q}}^{1/q}\leq\CP\sum_{k=-N}^N\PAR{\tau_k-\tau_{k+1}}^{2}\;\Capa_{\nu}(\Omega_k,\Omega_{k+1})\;.
$$
Next consider an appropriate choice of $(\tau_k)_{k=-N}^N$ : for $k\in\BRA{-N,\ldots N}$, let
$$
\tau_k=\sum_{j=k}^{N}\PAR{\frac{\mu(\Omega_j)}{\Capa_\nu(\Omega_j,\Omega_{j+1})}}^{\frac 1{2\,(1-q)}}\;.
$$
We observe that $\tau_k-\tau_{k+1}=\PAR{\frac{\mu(\Omega_k)} {\Capa_\nu(\Omega_k,\Omega_{k+1})}}^{1/(2(1-q))}$ and
$$
\BRA{\,\sum_{k=-N}^{N-1}\frac{\mu(\Omega_k)^{1/(1-q)}}{\Capa_\nu(\Omega_k,\Omega_{k+1})^{q/(1-q)}}}^{(1-q)/q}\leq 2^{1/q}\,\CP\left(1+\mathcal R_N\right)
$$
with $\mathcal R_N:=\frac{\mu(\Omega_N)^{1/(1-q)}}{\Capa_\nu(\Omega_N,\Omega_{N+1})^{q/(1-q)}}\,\BRA{\,\sum_{k=-N}^{N-1}\frac{\mu(\Omega_k)^{1/(1-q)}}{\Capa_\nu(\Omega_k,\Omega_{k+1})^{q/(1-q)}}}^{-1}$. By taking the limit as $N$ goes to infinity, we obtain $\betaP\leq 2^{1/q}\,\CP$.

\medskip {\sl Proof of \/}(ii). Let $f$ be a smooth non-negative function on $M$ and take $q\in(0,1]$. For all $a\geq0$,
$$
\var{\mu}{f^q}\leq\int\PAR{f^q-a^q}^2\,d\mu\leq\int |f-a|^{2q}\,d\mu\;.
$$
With $a:=m(f)$, a median of $f$ with respect to $\mu$, define $F_+=(f-a)_+$ and $F_-=(f-a)_-=F_+-(f-a)$, so that
$$
\var{\mu}{f^q}\leq\int\PAR{f^q-a^q}^2\,d\mu\leq\int F_+^{2q}\,d\mu+\int F_-^{2q}\,d\mu\;.
$$
We recall that $m=m(f)$ is a median of $f$ with respect to the measure $\mu$ if and only if $\mu(\{f\geq m\})\geq 1/2$ and $\mu(\{f\leq m\})\geq 1/2$. The computation of the term $\int F_-^{2q}\,d\mu$ is exactly the same as the one of $\int F_+^{2q}\,d\mu$, so we shall only detail one of them. Let us fix $\rho\in(0,1)$, note $\Omega_k:=\BRA{F_+\geq\rho^k}$ for any $k\in\dZ$, and use again the co-area formula:
\[
\int F_+^{2q}\,d\mu=\int_0^{+\infty}\kern -3pt\mu\PAR{\{F_+\geq t\}}\,d(t^{2q}) = \sum_{k\in\dZ}\int_{\rho^{k+1}}^{\rho^k}\mu\PAR{\{F_+\geq t\}}\,d(t^{2q}) \leq \frac{1-\rho^{2q}}{\rho^{2q}}\sum_{k\in\dZ}\mu\PAR{\Omega_k}\,\rho^{2kq}\,.
\]
By H\"older's inequality with parameters $(1/(1-q),1/q)$ one gets
\begin{equation*}
\begin{array}{rl}
\disp \int F_+^{2q}\,d\mu&\disp \leq\frac{1-\rho^{2q}}{\rho^{2q}}\,\PAR{\;\sum_{k\in\dZ}\frac{\mu(\Omega_k)^{1/(1-q)}}
{\Capa_\nu(\Omega_k,\Omega_{k+1})^{q/(1-q)}}}^{1-q}{\PAR{\;\sum_{k\in\dZ}\rho^{2k}\,
\Capa_\nu(\Omega_k,\Omega_{k+1})}^q}\\
&\disp \leq\frac{1-\rho^{2q}}{\rho^{2q}}\,\betaP^q {\PAR{\;\sum_{k\in\dZ}\rho^{2k}
\Capa_\nu(\Omega_k,\Omega_{k+1})}^q}.
\end{array}
\end{equation*}
For $k\in\dZ$, define $g_k:=\min\BRA{1,\big({\frac{F_+-{\rho^{k+1}}}{{\rho^k}-{\rho^{k+1}}}}\big)_+}$. Then we have $\1_{\Omega_k}\leq g_k\leq\1_{\Omega_{k+1}}$,
$$
\Capa_{\nu}(\Omega_k,\Omega_{k+1})\leq\int_{\Omega_{k+1}\setminus\Omega_k}\ABS{\nabla g_k}^2\,d\nu=\frac{1}{\rho^{2k}(1-\rho)^2}\int_{\Omega_{k+1}\setminus\Omega_k}\ABS{\nabla F_+}^2\,d\nu\;.
$$
Hence
$$
\int F_+^{2q}\,d\mu\leq\frac{1-\rho^{2q}}{\rho^{2q}\,(1-\rho)^{2q}}\,\betaP^q \PAR{\int\ABS{\nabla F_+}^2\,d\nu}^q\;.
$$
The same inequality holds for $F_-$:
$$
\int F _-^{2q}\,d\mu\leq\frac{1-\rho^{2q}}{\rho^{2q}\,(1-\rho)^{2q}}\,\betaP^q \PAR{\int\ABS{\nabla F _-}^2\,d\nu}^q\;.
$$
Using the inequality $a^q+b^q\leq 2^{1-q}(a+b)^q$ for any $a$, $b\geq0$, ones gets
$$
\PAR{\var{\mu}{f^q}}^{1/q}\leq \kappaP\,\betaP\int\ABS{\nabla f}^2\,d\nu
$$
with $\kappaP:=2^{(1-q)/q}\,\min_{\rho\in(0,1)}\frac{\PAR{1-\rho^{2q}}^{1/q}}{\rho^2\,(1-\rho)^2}$.
\end{eproof}

\subsection{$L^q$-logarithmic Sobolev inequalities}
\label{sec-logsob}

\begin{edefi}
Let $\mu$ and $\nu$ be respectively a probability measure and a positive measure on~$M$ and assume that $q\in(0,1]$. We shall say that $(\mu,\nu)$ satisfies a {\rm $L^q$-logarithmic Sobolev inequality\/} with constant $\CLSI$ if and only if, for any non-negative function $f\in{\mathcal C}^1(M)$,
\begin{equation*}
\label{eq2}
\ent{\mu}{f^{2q}}^{1/q}:=\PAR{\int f^{2q}\,\frac{\log f^{2q}}{\int f^{2q}\,d\mu}\,d\mu}^{1/q}\leq \CLSI \int\ABS{\nabla f}^2\,d\nu\;.
\end{equation*}
\end{edefi}
It is well known that $\ent{\mu}{f^2}\geq\var{\mu}f$ for any non-negative function $f$, for any probability measure $\mu$. Hence, for any $q\in(0,1]$, any $L^q$-logarithmic Sobolev inequality results in a $L^q$-Poincar\'e inequality with corresponding measures.

\medskip Let $q\in(0,1)$ and define
\begin{equation*}
\label{eq-beta2}
\betaSL=\sup\BRA{\sum_{k\in\dZ}\frac{\SBRA{\mu(\Omega_k)\,\log\!\PAR{1+\frac{e^2}{\mu(\Omega_k)}}}^{1/(1-q)}}{\SBRA{\Capa_\nu(\Omega_k,\Omega_{k+1})}^{q/(1-q)}}}^{(1-q)/q}
\end{equation*}
where the supremum is taken over all $\Omega\subset M$ with $\mu(\Omega)\le 1/2$ and all sequence
$\PAR{\Omega_k}_{k\in\dZ}$ such that, for all $k\in\dZ$, $\Omega_k\subset\Omega_{k+1}\subset \Omega$.

\begin{ethm}
\label{thm-ls}
Let $\mu$ and $\nu$ be respectively a probability measure and a positive measure on~$M$. If $q\in(0,1)$ and $\betaSL<+\infty$, then $(\mu,\nu)$ satisfies a $L^q$-logarithmic Sobolev inequality with constant $\CLSI\leq \kappaLS\,\betaSL$, where $\kappaLS$ depends only on $q$.
\end{ethm}
This theorem is the counterpart for the $L^q$-logarithmic Sobolev inequality of Theorem~\ref{thm-1},~(ii). As for Theorem~\ref{thm-1}, (i), related results will be stated in Corollary~\ref{cor-10}.

\medskip
\begin{eproof}
Let $f$ be a smooth function on $M$, $m=m(f)$ a median of $f$ with respect to $\mu$, and $\Omega_+:=\BRA{\ABS{f}>m}$, $\Omega_-:=\BRA{\ABS{f}<m}$. As in \cite{barthe-roberto}, we can write the dual formulation
\begin{multline}
\label{eq-bbr}
\ent{\mu}{f^{2q}}\leq\sup\BRA{\int (\ABS{f}^q-m^q)_+^{2}\,h\,d\mu\;:\;h\geq0\,,\;\int e^h\,d\mu\leq e^2+1}\\
+\sup\BRA{\int (\ABS{f}^q-m^q)_-^{2}\,h\,d\mu\;:\;h\geq0\,,\;\int e^h\,d\mu\leq e^2+1}\;.
\end{multline}
Such an inequality follows from Rothaus' estimate, \cite{Rothaus85},
$$
\ent{\mu}{g^2}\leq \ent{\mu}{(g-a)^2}+2\,\mu\((g-a)^2\)
$$
for any $a\in\R$, and the fact that, according to Lemma~5 in \cite{barthe-roberto},
$$
\ent{\mu}{(g-a)_+^2}+2\,\mu((g-a)_+^2)\leq\sup\BRA{\int (g-a)_+^{2}\,h\,d\mu\;:\;h\geq0\,,\;\int e^h\,d\mu\leq e^2+1}\;.
$$
Estimates for the positive and the negative part are exactly the same, so we will give details only for $F_+:=(\ABS{f}-m)_+$. Using the fact that $(t^q-1)^2<(t-1)^{2q}$ for any $t>1$, for $q\in(0,1)$, we get
$$
\int (\ABS{f}^q-m^q)_+^{2}\,h\,d\mu\leq\int F_+^{2q}\,h\,d\mu\;.
$$
Let $\rho\in(0,1)$.
\begin{multline*}
 \int F_+^{2q}\,h\,d\mu=\int_0^{+\infty}\int_{F_+>t}\,h\,d\mu\,d(t^{2q})=\sum_{k\in\dZ}\int_{\rho^{k+1}}^{\rho^k}\int_{F_+>t}\,h\,d\mu\,d(t^{2q})\\
\leq\sum_{k\in\dZ}\PAR{\rho^{2qk}-\rho^{2q(k+1)}}\int_{F_+>\rho^{k+1}}\,h\,d\mu=\frac{1-\rho^{2q}}{\rho^{2q}}\sum_{k\in\dZ}\rho^{2qk}\int_{F_+>\rho^k}\,h\,d\mu\;.
\end{multline*}
Using Lemma~6 of \cite{barthe-roberto}, which asserts that
$$
\int_{F_+>\rho^k}\,h\,d\mu\leq\mu(\Omega_k)\,\log\!\PAR{1+\frac{e^2}{\mu(\Omega_k)}}
$$
where $\Omega_k:=\BRA{F_+>\rho^k}$, we obtain
\begin{equation*}
\sup\BRA{\int F_+^{2q}\,h\,d\mu\;:\;h\geq0\,,\;\int e^h\,d\mu\leq e^2+1}\leq\frac{1-\rho^{2q}}{\rho^{2q}} \sum_{k\in\dZ}\rho^{2qk}\mu(\Omega_k)\,\log\!\PAR{1+\frac{e^2}{\mu(\Omega_k)}}\,.
\end{equation*}
By H\"older's inequality, it follows that
\begin{equation*}
\begin{array}{rl}
\disp\ent{\mu}{F_+^{2q}}\!\!&\!\disp \leq\frac{1-\rho^{2q}}{\rho^{2q}}\sum_{k\in\dZ}\rho^{2qk}\,\mu(\Omega_k)\,\log\!\PAR{1+\frac{e^2}{\mu(\Omega_k)}}\\
&\!\disp \leq\!\frac{1-\rho^{2q}}{\rho^{2q}}\,\PAR{\sum_{k\in\dZ} \frac{\SBRA{\mu(\Omega_k)\,\log\!\PAR{1+\frac{e^2}{\mu(\Omega_k)}}}^{1/(1-q)}} {\SBRA{\Capa_\nu(\Omega_k,\Omega_{k+1})}^{q/(1-q)}}}^{\kern -3pt 1-q}{\PAR{\,\sum_{k\in\dZ}\rho^{2k}\,\Capa_\nu(\Omega_k,\Omega_{k+1})}^{\kern -3pt q}}\\
&\!\disp \leq\frac{1-\rho^{2q}}{\rho^{2q}}\,\betaSL^q {\PAR{\,\sum_{k\in\dZ}\rho^{2k}\,\Capa_\nu(\Omega_k,\Omega_{k+1})}^{\kern -3pt q}}\\
&\!\disp\leq\frac{1-\rho^{2q}}{\rho^{2q}}\,\betaSL^q\PAR{\int\ABS{\nabla F_+}^2\,d\mu}^{\kern -3pt q}\,.
\end{array}
\end{equation*}
The same computation shows that
$$
\ent{\mu}{F_-^{2q}}\leq\frac{1-\rho^{2q}}{\rho^{2q}}\,\betaSL^q\PAR{\int\ABS{\nabla F_-}^2\,d\mu}^{\kern -3pt q}\,.
$$
Summing both contributions in Inequality~\eqref{eq-bbr} completes the proof with
$$
\kappaLS=2^\frac{1-q}q\PAR{1-\rho^{2q}}^{1/q}\rho^{-2}\;.
$$
\end{eproof}

\section{Weak inequalities and explicit criteria}
\label{Further}

The goal of this section is to provide tractable criteria to establish $L^q$-Poincar\'e and the $L^q$-logarithmic Sobolev inequalities. The strategy here is to adapt results which have been obtained for weak Poincar\'e inequalities by Barthe, Cattiaux and Roberto in \cite{ca-ba-ro2}. Two important results stated in this paper are extended to measures $\mu$ and $\nu$ which are not supposed to be absolutely continuous with respect to the volume measure, and given with proofs in Section~\ref{Sec:Appendix}.

\subsection{$L^q$ Poincar\'e and weak Poincar\'e inequalities}
\label{Sec:WP}

Even if the constants $\betaP$ and $\betaSL$ provide an estimate of the best constant of the $L^q$-Poincar\'e and the $L^q$-logarithmic Sobolev inequalities, their expressions in terms of suprema taken over infinite sequences of sets are {\it a priori\/} difficult to use. In this section, we look for simpler criteria and establish upper and lower bounds on the constants.

\medskip The first idea is relate the $L^q$-Poincar\'e inequality and the weak Poincar\'e inequality introduced by R\"ockner and the fourth author in~\cite{r-w}. Let us define the {\sl oscillation\/} of a bounded function $f$ by $\osc{f}:={\rm supess}_\mu f-{\rm infess}_\mu f$. If $\mu$ is absolutely continuous with respect to the volume measure and $f$ is continuous, we can therefore define such a quantity as $({\rm sup} \tilde f-{\rm inf} \tilde f)$ where $\tilde f$ is the restriction of $f$ to the support of~$\mu$. Our definition slightly differs from the one of \cite{r-w}, which is based on ${\rm supess}_\mu |f-\int f\,d\mu|$.
\begin{edefi}
\label{Defn:WeakPoincare}
Let $\mu$ and $\nu$ be respectively a probability measure and a positive measure on~$M$. We shall say that $(\mu,\nu)$ satisfies a {\rm weak Poincar\'e inequality} if there exists a non-negative non increasing function on $(0,+\infty)\!\ni\!s\mapsto\CWP(s)$ such that, for any bounded function $f\in{\mathcal C}^1(M)$,
$$
\forall\;s>0\,,\quad\var{\mu}f\leq\CWP\PAR{s}\int\ABS{\nabla f}^2 \,d\nu+s\,\big[\osc{f}\big]^2\,.
$$
\end{edefi}
Since $\var{\mu}f\leq \mu((f-a)^2)$ for all $a\in\R$, and as a special case, for $a=({\rm supess}_\mu f+{\rm infess}_\mu f)/2$, $\var{\mu}f\leq\big[\osc{f}\big]^2/4$, which means that we can assume that $\CWP\PAR{s}\equiv 0$ for any \hbox{$s\geq 1/4$}.
\begin{eprop}
\label{cor-1}
Let $q\in[1/2,1)$. If $(\mu,\nu)$ satisfies the $L^q$-Poincar\'e inequality, then it also satisfies a weak Poincar\'e inequality with $\CWP(s)=\kkappa\,\betaP\,s^{1-1/q}$, $\kkappa:=(11+5\sqrt 5)/2$.
\end{eprop}
\begin{eproof}
By Theorem~\ref{thm-1}, the constant $\betaP$ is finite. Let $A\subset \Omega\subset M$ with $\mu(\Omega)\leq1/2$ and consider the sequence: $\PAR{\Omega_k}_{k\in\dZ}$ such that $\Omega_k=A$ for all $k\leq0$ $\Omega_k=\Omega$ and for all $k>0$. Notice that $\Capa_{\nu}(\Omega_k,\Omega_{k+1})=\infty$ if $\Omega_k=\Omega_{k+1}$. By definition of $\betaP$ we get
$$
\frac{\SBRA{\mu(A)}^{1/q}}{\betaP}\leq\Capa_{\nu}(A,\Omega)\;.
$$
Using the method of Barthe, Cattiaux and Roberto in \cite{ca-ba-ro2}, Theorem~2, one can then prove that $(\mu,\nu)$ satisfies a weak Poincar\'e inequality with constant $\CWP(s)$. See Theorem~\ref{crit-capa} for a precise statement and apply it with $\gamma(s)=\betaP\,s^{1-1/q}$.
\end{eproof}

\medskip Another criterion to prove $L^q$-Poincar\'e inequalities is based on Theorem~2.3.6 of \cite{mazja}.
\begin{ethm}[\cite{mazja}]
\label{thm-maz}
Let $q\in[1/2,1)$. For all bounded open set $\Omega\subset M$, if $\PAR{\Omega_k}_{k\in\dZ}$ is an increasing sequence of open sets such that $\Omega_k\subset\Omega_{k+1}\subset \Omega$, then
$$
{\sum_{k\in\dZ}\frac{\mu\PAR{\Omega_k}^{1/(1-q)}}{\SBRA{\Capa_\nu(\Omega_k,\Omega_{k+1})}^{q/(1-q)}}}\leq \frac{1}{1-q}\int_0^{\mu(\Omega)}\PAR{\frac{t}{\Phi(t)}}^{q/(1-q)}\,dt\;,
$$
where $\Phi(t):=\inf\BRA{\Capa_\nu(A,\Omega)\,:\,A\subset \Omega\,,\;\mu(A)\geq t}$.
\end{ethm}
Notice that as a consequence, $\betaP\leq (1-q)^{-(1-q)/q}\,\|t/\Phi(t)\|_{L^{q/(1-q)}(0,\mu(\Omega))}$
\begin{ecor}
\label{cor-2}
Let $q\in[1/2,1)$ and assume that $(\mu,\nu)$ satisfies a weak Poincar\'e inequality with function $\CWP$. Then $(\mu,\nu)$ satisfies satisfies a $L^q$-Poincar\'e inequality with
$$
\betaP\leq\kappaP\,\PAR{\scriptstyle \frac{4}{1-q}}^{\frac{1-q}q}\,\|\CWP(\cdot/4)\|_{L^{\frac q{1-q}}(0,1/2)}
$$
where $\kappaP$ is defined in Theorem~\ref{thm-1}.\end{ecor}
\begin{eproof}
The method is again similar to the one of Theorem~2 in \cite{ca-ba-ro2}; see Theorem~\ref{crit-capa} in Section~\ref{Sec:Appendix}. If $(\mu,\nu)$ satisfies a weak Poincar\'e inequality, then for all $\Omega\subset M$ with $\mu(\Omega)\le 1/2$, $A\subset \Omega$,
$$
\frac{\mu(A)}{4\,\CWP\PAR{\mu(A)/4}}\leq\Capa_\nu(A,\Omega)\;.
$$
Hence for $t>0$, $\Phi(t)\geq \frac{t}{4\,\CWP\PAR{t/4}}$, and the result follows.
\end{eproof}

\medskip Proposition~\ref{cor-1} and Corollary~\ref{cor-2} can be summarized as follows. For any $q\in[1/2,1)$,
$$
L^q\mbox{-Poincar{\'e}} \quad\Longrightarrow\quad \begin{array}{c}\mbox{Weak Poincar\'e}\\ \mbox{with }\CWP(s)=C\,s^{\frac{q-1}q}\end{array} \quad\Longrightarrow\quad \begin{array}{c}L^{q'}\mbox{-Poincar{\'e}}\\ \forall\;q'\in(0,q)\end{array}\;.
$$
As we shall see in Section~\ref{sec-link}, weak Poincar\'e inequalities with $\CWP(s)=C\,s^{\frac{q-1}q}$ do not imply $L^{q'}$-Poincar\'e inequalities with $q'=q$.

\subsection{$L^q$-Logarithmic Sobolev and weak logarithmic Sobolev inequalities}
\label{Sec:WLS}

$L^q$-Poincar\'e have been established in terms of weak Poincar\'e inequalities in Section~\ref{Sec:WP}. Very similar characterizations can be done for $L^q$-logarithmic Sobolev in terms of weak logarithmic Sobolev inequalities. Recall first the definition of the weak logarithmic Sobolev inequality.
\begin{edefi}
Let $\mu$ and $\nu$ be respectively a probability measure and a positive measure on~$M$. We sall say that $(\mu,\nu)$ satisfies a {\rm weak logarithmic Sobolev inequality} if there exists a positive and non-increasing function $\CWL$ on $\dR^+$ such that for any bounded function $f\in{\mathcal C}^1(M)$,
$$
\forall\;s>0\,,\quad\ent{\mu}{f^2}\leq\CWL\PAR{s}\int\ABS{\nabla f}^2 \,d\nu+s\,\big[\osc{f}\big]^2\;.
$$
\end{edefi}
A preliminary step amounts to state the analogue of Proposition \ref{cor-1}.
\begin{eprop}
\label{cor-1ls}
Let $q\in[1/2,1)$. If $\mu$ is absolutely continuous with respect to the volume measure and $(\mu,\nu)$ satisfies the $L^q$-logarithmic Sobolev inequality, then it also satisfies a weak logarithmic Sobolev inequality with $\CWL(s)= c_q\, s^{1-1/q}$ for some positive constant $c_q$.
\end{eprop}
\begin{eproof}
By Legendre duality, for any non-negative function $f$,
$$
\int f^{2q}\,g\;d\mu\leq\ent{\mu}{f^{2q}}\quad\forall\;g\;\mbox{such that}\;\int e^g\,d\mu\leq 1\;.
$$
Let $A\subset \Omega\subset M$ with $\mu(\Omega)\leq 1$ and assume that $f\in{\mathcal C}^1(M)$ is such that $\1_A\leq f\leq\1_\Omega$. Then by the $L^q$-logarithmic Sobolev inequality for $(\mu,\nu)$, we get
$$
\PAR{\int_Ag\;d\mu}^{1/q}\le\PAR{\int f^{2q}\,g\;d\mu}^{1/q}\le\betaSL\int\ABS{\nabla f}^2\,d\nu\;.
$$
Choose now $g=-\infty$ on $\Omega^c$, $g=0$ on $\Omega\backslash A$ and $g=\log(1+1/(2\mu(A))$ on $A$ so that $\int e^gd\mu\le 1$. Using $\1_A\le f\le \1_{\Omega}$, a simple computation gives
$$
\ent{\mu}{f^{2q}}^{1/q}\ge \SBRA{\mu(A)\log\left(1+\frac{1}{2\,\mu(A)}\right)}^{1/q}\ge c^{1/q}\frac{\mu(A)\log\left(1+\frac{e^2}{\mu(A)}\right)}{\SBRA{\mu(A)\log\left(1+\frac{e^2}{\mu(A)}\right)}^{1-1/q}}
$$
where $c:=\log 2/\log(1+2e^2)$. Now optimizing in $f$ leads to
$$
\Capa_\nu(A,\Omega)\ge \frac{c^{1/q}}{\betaSL}\frac{\mu(A)\log\left(1+\frac{e^2}{\mu(A)}\right)}{\SBRA{\mu(A)\log\left(1+\frac{e^2}{\mu(A)}\right)}^{1-1/q}}\;,
$$
which entails the desired weak logarithmic Sobolev inequality by \cite[Th. 2.2]{ca-ge-gu}
\end{eproof}

\medskip
Let us study the reciprocal property.
\begin{ethm}
\label{prop2}
Let $q\in[1/2,1)$. For any $\Omega\subset M$, $\PAR{\Omega_k}_{k\in\dZ}$ such that $\Omega_k\subset\Omega_{k+1}\subset \Omega$, one gets
$$
{\sum_{k\in\dZ}\frac{\SBRA{\mu\PAR{\Omega_k}\log\!\PAR{1+\frac{e^2}{\mu\PAR{\Omega_k}}}}^{1/(1-q)}} {\SBRA{\Capa_\nu(\Omega_k,\Omega_{k+1})}^{q/(1-q)}}}\leq\frac{1}{1-q}\int_0^{\mu(\Omega)}\PAR{\frac{t}{\psi(t)}}^{q/(1-q)}\,dt\;,
$$
where $\psi(t):=\inf\BRA{\Capa_\nu(A,\Omega)\,:\,A\subset \Omega\,,\;\mu(A)\,\log\!\PAR{1+\frac{e^2}{\mu(A)}}\geq t}$.
\end{ethm}
\begin{eproof}
The proof is a direct adaptation of the proof of Theorem~\ref{thm-maz} (see pages~122 and 123 of \cite{mazja}).
\end{eproof}

\medskip
As a consequence, we obtain the following characterization. Let $c:=\log 2/(2\log(1+2e^2))$.
\begin{ecor}
\label{cor-10}
Let $q\in[1/2,1)$. Assume that $(\mu,\nu)$ satisfies a weak logarithmic Sobolev inequality with function $\CWL$. Then it satisfies a $L^q$-logarithmic Sobolev inequality~if
$$
\int_0^{1/2}{{\CWL(c\,t)}}^{q/(1-q)}\,dt<+\infty\;.
$$
\end{ecor}
In such a case, the optimal constant of the $L^q$-logarithmic Sobolev inequality is bounded by $\kappaLS\big({\frac{4}{1-q}\int_0^{1/2}{{\CWL(c\,t)}}^{q/(1-q)}\,dt}\big)^{(1-q)/q}$, where $\kappaLS$ is defined in Theorem~\ref{thm-ls}.

\medskip \begin{eproof}
By Theorem~2.2 in \cite{ca-ge-gu} (also see~Lemma~\ref{Lem:CGG05-2.2}), for all $\Omega\subset M$ with $\mu(\Omega)\le 1/2$, $A\subset \Omega$,
$$
 \frac{\frac{\mu(A)}{2}\log\!\PAR{1+\frac{1}{2\mu(A)}}}{\CWL\SBRA{\frac{\mu(A)}{2}\log\!\PAR{1+\frac{1}{2\mu(A)}}}} \leq\Capa_\nu(A,\Omega)\;,
$$
Picking $c=\log(2)/(2\log(1+2e^2))$, we thus get
$$
c\,\frac{{\mu(A)}\log\!\PAR{1+\frac{e^2}{\mu(A)}}}{\CWL\SBRA{c\,{\mu(A)}\log\!\PAR{1+\frac{e^2}{\mu(A)}}}}\leq \Capa_\nu(A,\Omega)\;.
$$
Theorem~\ref{prop2} then gives the result because $\psi(t)\ge c\,t\,/\,\CWL(c\,t)$.
\end{eproof}

\medskip 
As a consequence of Proposition \ref{cor-1ls} and Corollary \ref{cor-10}, we have the following result.
\begin{ecor}
\label{cor-10bis}
Let $q\in[1/2,1)$ and assume that $(\mu,\nu)$ satisfies a $L^q$-logarithmic Sobolev inequality. Then for all $0<q'<q$, all $\Omega\subset M$ with $\mu(\Omega)\le 1/2$ and all sequence $\PAR{\Omega_k}_{k\in\dZ}$ such that, for all $k\in\dZ$, $\Omega_k\subset\Omega_{k+1}\subset \Omega$, we have
\begin{equation*}
\betaSL^{q'}=\sup\BRA{\sum_{k\in\dZ}\frac{\SBRA{\mu(\Omega_k)\,\log\!\PAR{1+\frac{e^2}{\mu(\Omega_k)}}}^{1/(1-q')}} {\SBRA{\Capa_\nu(\Omega_k,\Omega_{k+1})}^{q'/(1-q')}}}^{(1-q')/q'}<+\infty\;.
\end{equation*}
\end{ecor}
This result completes that of Theorem~\ref{thm-ls}. Unfortunately the equivalence is not proved for the $L^q$-logarithmic Sobolev inequality. This however proves the counterpart of Proposition~\ref{prop-hi}, namely the hierarchy between $L^q$-logarithmic Sobolev inequalities. Summarizing the results of this subsection, we have for any $q\in[1/2,1)$,
$$
L^q\mbox{-Logarithmic Sobolev} \quad\!\Longrightarrow \begin{array}{c}\mbox{Weak logarithmic Sobolev}\\ \mbox{with }\CWL(s)=C\,s^{\frac{q-1}q}\end{array}\Longrightarrow \begin{array}{c}L^{q'}\mbox{-Logarithmic Sobolev}\\ \forall\;q'\in(0,q)\end{array}\;.
$$

\subsection{A Hardy condition on $\dR$}
\label{sec-dim1}

On $\dR$, to a probability measure $\mu$ and a positive measure $\nu$ with density $\rho_\nu$ with respect to Lebesgue's measure, if $m_\mu$ is a median of $\mu$, we associate the functions
$$
R(x):=\mu([x,+\infty))\;,\quad L(x):=\mu((-\infty,x])\;,\quad r(x):=\!\int_{m_\mu}^x\frac 1{\rho_\nu}\;dx\quad\mbox{and}\quad \ell(x):=\!\int_x^{m_\mu}\!\frac 1{\rho_\nu}\;dx\;.
$$
\begin{eprop}
\label{thm-dim1}
Let $q\in[1/2,1]$, and let $\mu$ and $\nu$ be respectively a probability measure and a positive measure on~$\dR$. With the above notations, $(\mu,\nu)$ satisfies a $L^q$-Poincar\'e inequality if
$$
\int_{m_\mu}^\infty |\,r\,R\,|^{q/(1-q)}\,d\mu<\infty\quad\mbox{and}\quad\int_{-\infty}^{m_\mu} |\,\ell\,L\,|^{q/(1-q)}\,d\mu<\infty\;.
$$
Analogoulsy, $(\mu,\nu)$ satisfies a $L^q$-logarithmic Sobolev inequality if
$$
\int_{m_\mu}^\infty |\,r\,R\,\log R\,|^{q/(1-q)}\;d\mu<\infty\quad\mbox{and}\quad\int_{-\infty}^{m_\mu} |\,\ell\,L\,\log L\,|^{q/(1-q)}\;d\mu<\infty\;.
$$
\end{eprop}
\begin{eproof}
The proof of Theorem~3 in \cite{ca-ba-ro2} can then be adapted to the setting of Proposition~\ref{thm-dim1}. It relies on weak Poincar\'e inequalities. Taking advantage of $\var\mu f\leq \mu(|F_-|^2)+\mu(|F_+|^2))$ with $F_\pm:=(f-f(m_\mu))_\pm$, we notice that the weak Poincar\'e inequality 
$$
\var\mu f\leq\kkappa\,\gamma(s)\int\ABS{\nabla f}^2 \,d\nu+s\,\big[\osc{f}\big]^2\quad\forall\;s\in (0,1/2)
$$
holds if we are able to prove independently for $g=F_+$ and $g=F_-$ that the inequality
$$
\mu(|g|^2)\leq\kkappa\,\gamma(s)\int\ABS{\nabla g}^2 \,d\nu+s\,\big[\mbox{\rm supess}_\mu{g}\big]^2\quad\forall\;s\in (0,1/2)
$$
holds for some positive non increasing function $\gamma$ on $(0,1)$ and for $\kkappa:=(11+5\sqrt 5)/2$. For this purpose, we are going to rely on Lemma~\ref{prel-crit-capa}. If $A$ and $B$ are two measurable subsets of $M=(m_\mu,\infty)$ such that $A\subset B$ and $\mu(B)\leq 1/2$, then
$$
\Capa_\nu(A,B)\geq\Capa_\nu\big(A,(m_\mu,\infty)\big)=\Capa_\nu\big((a,\infty),(m_\mu,\infty)\big)=\frac 1{r(a)}
$$
where $a=\inf A$. By Lemma~\ref{prel-crit-capa}, it is therefore sufficient to prove that
$$
\frac 1{r(a)}\geq\frac{R(a)}{\gamma(R(a))}\quad\forall\;a>m_\mu\;.
$$
With the change of variables $t=R(a)$, $a>m_\mu$, this amounts to require that
$$
\gamma(t)\geq t\, (r\circ R^{-1})(t)\;.
$$
With no restriction, we can choose $\gamma(t):=t\,(r\circ R)^{-1}(t)$ for any $t\in(0,1/2)$. 

\medskip
By Corollaries~\ref{cor-2} and \ref{cor-10}, $(\mu,\nu)$ satisfies a $L^q$-Poincar\'e inequality if $\CWP\in L^{q/(1-q)}(0,1/2)$ and a $L^q$-logarithmic inequality if $\CWL\in L^{q/(1-q)}(0,1/2)$. 
\end{eproof}

\subsection{Examples}

Let us illustrate the above results on $L^q$-Poincar\'e and $L^q$-logarithmic Sobolev inequalities with examples on $M=\dR$, in case of a single measure $\mu=\nu$. We start with some observations on $L^q$-Poincar\'e inequalities.
\begin{enumerate}
\item[{\rm (i)}] The classical Poincar\'e inequality implies a $L^q$-Poincar\'e inequality for all $q\in[1/2,1)$ by Proposition~\ref{prop-hi}. This gives an explicit estimate of the constant $\kappa_{p,q}$ of Theorem~1 of \cite{car-al} when $p\,q=2$, that is $\kappa_{p,q}\geq 2^{p+2}\,\pi^2$. Recall indeed that $1/(4\,\pi^2)$ is the Poincar\'e constant of the uniform measure on $[0,1)$, with periodic boundary conditions.
\item[{\rm (ii)}] For $p\in (0,1)$, the probability measure $d\mu=e^{-\ABS{x}^p}/(2\,\Gamma(1+1/p))\,dx$, $x\in\dR$, satisfies a weak Poincar\'e inequality with $\CWP(s)=C\log (2/s)^{2/p-2}$ for some positive constant $C$. As a consequence, $\mu$ also satisfies a $L^q$-Poincar\'e inequality for all $q\in[1/2,1)$.
\item[{\rm (iii)}] For $\al>0$, the probability measure $d\mu=\al\,(1+\ABS{x})^{-1-\al}\,dx/2$, $x\in\dR$, satisfies a weak Poincar\'e inequality with $\CWP(s)=C\,s^{-2/\al}$ for some positive constant $C$, see \cite{r-w,ca-ba-ro2}. Then for any $q\in[1/2,1)$, the probability measure $\mu$ satisfies a $L^q$-Poincar\'e inequality if $\al>2q/(1-q)$. As in Example~1.1 of \cite{2006-wang} in case of Orlicz-Poincar\'e inequalities (see below), the $L^q$-Poincar\'e inequality is not satisfied for $\al=2q/(1-q)$. 
\end{enumerate}

\medskip
Similar remarks can be done for $L^q$-logarithmic Sobolev inequalities.
\begin{enumerate}
\item [{\rm (i)}] Gross' logarithmic Sobolev inequality implies the $L^q$-logarithmic Sobolev inequality, for all $q\in[1/2,1)$.
\item [{\rm (ii)}] For $p\in (0,1)$, the probability measure $d\mu=e^{-\ABS{x}^p}/(2\,\Gamma(1+1/p))\,dx$, $x\in\dR$, satisfies a weak logarithmic Sobolev inequality with $\CWL(s)=C\,\PAR{\log 1/s}^{(2-p)/p}$ for some positive constant $C$, see \cite{ca-ge-gu}. As a consequence, $d\mu$ also satisfies a $L^q$-logarithmic Sobolev inequality for all $q\in[1/2,1)$.
\item [{\rm (iii)}] For $\al>0$, the probability measure $d\mu=\al\,(1+\ABS{x})^{-1-\al}\,dx/2$, $x\in\dR$, satisfies the weak logarithmic Sobolev with $\CWL(s)=C\,s^{-2/\al}\,\PAR{\log(1/s)}^{(2+\al)/\al}$ for some positive constant $C$, see \cite{ca-ge-gu}. Then for any $q\in[1/2,1)$, the probability measure $\mu$ satisfies a $L^q$-logarithmic Sobolev inequality if $\al>2q/(1-q)$.
\end{enumerate}

\medskip
At the light of the above examples $L^q$-Poincar\'e and $L^q$-logarithmic Sobolev inequalities seem to be satisfied by the same measures. This is not true as shown by the following example. On~$\dR$, the probability measure
$$
d\mu=\frac{C_{\al,\be}}{1+\ABS{x}^{1+\al}\,\ABS{\log x}^{\be}}\;dx\quad\mbox{with}\quad\al>0\,,\;\beta\in\dR\;,
$$
satisfies a weak Poincar\'e inequality with $\CWP(s)\!=\!C\,s^{-2/\al}\PAR{\log(1/s)}^{-2\be/\al}$ for some constant $C>0$ and a weak logarithmic Sobolev inequality with $\CWL(s)\!=\!C's^{-2/\al}\PAR{\log(1/s)}^{1+2(1-\be)/\al}$ for some positive constant $C'$. Fix $\al$ such that $\frac{2}{\al}\frac{q}{1-q}=1$. Using Bertrand's integrals, the probability measure $\mu$ satisfies a $L^q$-Poincar\'e inequality if and only if $\be>1$, and a $L^q$-logarithmic Sobolev inequality if and only if $\be>1+1/(1-q)$. The two conditions clearly differ.

\subsection{Orlicz-Poincar\'e inequalities}
\label{sec-link}

The $L^q$-Poincar\'e inequality for $q\in[1/2,1)$ is a particular case of the Orlicz-Poincar\'e inequality introduced by Roberto and Zegarlinski in \cite{ro-zeg} in the sub-Gaussian case and by the fourth author in \cite{2006-wang}, in the others cases.
\begin{eprop} \label{prop-wang}
Let $q\in [1/2, 1]$. Then the $L^q$-Poincar\'e inequality holds for some $\CP>0$ if and only if the following Orlicz-Poincar\'e inequality
\begin{equation}
\label{2}\PAR{\int\ABS{f-\mu(f)}^{2q}\,d\mu}^{1/q}\le\COP\int\ABS{\nabla f}^2\,d\mu\;.
\end{equation}
holds for some $\COP>0.$
\end{eprop}
\begin{eproof} We notice that $\var{\mu}{f^q}\le\int(f^q-a)^2\,d\mu$ for any $a\in\R$, and as a special case for $a=\mu(f)$. The function $t\mapsto(t^q-1)/(t-1)^q$ is monotone increasing on $(1,\infty)$ and converges to $1$ as $t\to\infty$, so that $(t^q-1)^2\leq(t-1)^{2q}$ for any $t\in(1,\infty)$. This proves that $\var{\mu}{f^q}\le \int |f-\mu(f)|^{2q}\,d\mu$. Inequality~(\ref{2}) therefore implies the $L^q$-Poincar\'e inequality. On the other hand, let $F:= f-m$, where $m$ is a median of $f$. We have
$$
\|f-\mu(f)\|_{2q} = \|F-\mu(F)\|_{2q}\le \|F\|_{2q}+\mu(|F|)\le 2\,\|F\|_{2q}\;,
$$
and hence
$$
\mu(|f-\mu(f)|^{2q})\le 2^{2q} \int(F_+^{2q}+F_-^{2q})\;d\mu\;.
$$
As in the proof of Theorem~\ref{thm-1}, (\ref{2}) follows from $\betaP<\infty$.
\end{eproof}

\medskip We are now in position to prove that weak Poincar\'e inequalities with $\CWP(s)=C\,s^{\frac{q-1}q}$ and $L^q$ Poincar\'e inequalities are not equivalent (see the end of Section~\ref{Sec:WP}), or to be precise that
$$
\begin{array}{c}\mbox{Weak Poincar\'e}\\ \mbox{with }\CWP(s)=C\,s^{\frac{q-1}q}\end{array} \quad\kern10pt\not{\kern-10pt}\Longrightarrow\quad L^q\mbox{-Poincar{\'e}}\;.
$$
By Proposition~\ref{prop-wang} and according to \cite[Proposition~3.2]{2006-wang}, $\mu(f^{2q/(1-q)})$ is finite for any $f\in\Cn{1}$ with $\|f\|_{{\rm Lip}}\le 1$. An example for which a weak Poincar\'e inequality with $\CWP(s)=C\,s^{(q-1)/q}$ holds while the $L^q$ Poincar\'e inequality is wrong is given by $\mu=\nu=\al\,(1+\ABS{x})^{-1-\al}\,dx/2$, $x\in\dR$, $\al=2q/(1-q)$, which satisfies a weak Poincar\'e inequality with $\CWP(s)=C\,s^{(q-1)/q}$ for some positive constant~$C$, and $f(x):=\sqrt{1+x^2}$.

\subsection{Perturbation, tensorization and concentration of measure}
\label{sec-ptc}

\begin{eprop}
\label{Prop:Pertub-Tensor-Concentr}
\begin{enumerate}[{\rm (i)}]
\item Let $\mu$, $\nu$ be respectively a probability measure and a positive measure on~$M$. Assume that $h$ is a bounded function on $M$ and define the probability measure $d\mu_h:=Z_h^{-1}\,e^h\,d\mu$ with $Z_h:=\int e^h\,d\mu$. If $(\mu,\nu)$ satisfies a $L^q$-Poincar\'e (resp. a $L^q$-logarithmic Sobolev) inequality with constant~$\CP$ (resp. $\CLSI$), then $(\mu_h,\nu)$ satisfies a $L^q$-Poincar\'e (resp. a $L^q$-logarithmic Sobolev) inequality with constant $\CP\exp\PAR{\osc{h}/q}$ (resp. $\CLSI\exp\PAR{\osc{h}/q}$).
\item If for any $i\in\BRA{1,\cdots, n}$, $\mu_i$ is a probability measure and $(\mu_i,\mu_i)$ satisfies a $L^q$-Poincar\'e (resp. a $L^q$-logarithmic Sobolev) inequality with constant ${\CP}_i$ (resp. ${\CLSI}_i$), then $(\otimes_{i=1}^n\mu_i$, $\otimes_{i=1}^n\mu_i)$ satisfies a $L^q$-Poincar\'e  (resp.  $L^q$-logarithmic Sobolev) inequality on $M^n$ with constant ${n}^{1/q-1}\!\max_{1\leq i\leq n}{\CP}_i$ (resp. ${n}^{1/q-1}\!\max_{1\leq i\leq n}{\CLSI}_i$).
\item  If $\mu$ is a probability measure and $(\mu,\mu)$ satisfies a $L^q$-Poincar\'e inequality, then for any non-negative function $f\in\Cn{1}$ with $\|f\|_{{\rm Lip}}\le 1$ there exists $t_0>0$ and $C>0$ such that
\begin{equation}
\label{eq-concen}
\forall\;t\geq t_0\,,\quad\mu(\BRA{f\ge t})\le \frac{C}{t^{2q/(1-q)}}\;.
\end{equation}
\end{enumerate}
\end{eprop}
\begin{eproof} 
The first point is based on the same proof as in Theorem~3.4.1 and~3.4.3 of \cite{logsob}. We observe that for any $a\in\dR$, 
$$
\var{\mu_h}{f^q}\leq\int\left|f^q-a\right|^2\,d\mu_h\leq e^{-\osc h}\int |f^q-a|^2\,d\mu\;,
$$
and apply the $L^q$-Poincar\'e inequality with $a=\mu(f)$. Similarly, with $a=\mu(f)$, we get
$$
\ent{\mu_h}{f^q}\leq\int\SBRA{f^q\,\log\PAR{\frac{f^q}a}+f^q-a}\;d\mu_h\leq e^{-\osc h}\,\ent{\mu}{f^q}\;.
$$

The second point is almost the same as in Theorem~3.2.1 and 3.2.2 of \cite{logsob}:  for all functions $f\in\Cn{1}(M)$ (see Proposition~1.4.1 of \cite{logsob}), by sub-additivity of the variance, 
$$
\var{\otimes_{i=1}^n\mu_i}{f^q}\leq\sum_{i=1}^n \int \var{\mu_i}{f^q}\,d\!\otimes_{i=1}^n\mu_i\;.
$$
Apply the $L^q$-Poincar\'e inequality, $\big[\var{\mu_i}{f^q}\big]^{1/q}\leq{\CP}_i\int \ABS{\nabla_{x_i} f}^2\,d\mu_i$, component by component to get
$$
\var{\otimes_{i=1}^n\mu_i}{f^q}\leq\sum_{i=1}^n \int\PAR{ {\CP}_i\int \ABS{\nabla_{x_i} f}^2\,d\mu_i}^q\,d\!\otimes_{i=1}^n\mu_i\;.
$$
H\"older's inequality with $q\in(0,1)$, and the identity $\sum_{i=1}^n x_i^q\leq n^{1-q}(\sum_{i=1}^n x_i)^q $ for $x_i\geq 0$, give 
$$
\var{\otimes_{i=1}^n\mu_i}{f^q}\leq\sum_{i=1}^n \PAR{{\CP}_i\int \ABS{\nabla f}^2\,d\!\otimes_{i=1}^n\mu_i}^q \leq n^{1-q}\PAR{\max_{1\leq i\leq n}{\CP}_i\,\sum_{i=1}^n \int{ \ABS{\nabla f}^2 }\,d\!\otimes_{i=1}^n\mu_i}^q\,,
$$
with the notation $\sum_{i=1}^n\ABS{\nabla_{x_i} f}^2=\ABS{\nabla f}^2$. The proof for the $L^q$-logarithmic Sobolev inequality is similar and relies on the sub-additivity of the entropy:
$$
\ent{\otimes_{i=1}^n\mu_i}{f^{2q}}\leq\sum_{i=1}^n \int \ent{\mu_i}{f^{2q}}\,d\!\otimes_{i=1}^n\mu_i\;.
$$

Property (iii) is inspired by the method of Aida, Masuda and Shigekawa in \cite{aida}. Define $a(t):=\mu(\{f\ge t\})$ and choose $t_0$ such that $a(t_0)\le 1/2$. For any $t\geq t_0$, define $g:=\min\BRA{\frac1t\,(f-t)_+,1}$. On the one hand, by the Cauchy-Schwarz inequality,
$$
\PAR{\int g^q\,d\mu}^2=\PAR{\int_{f\geq t_0}g^q\,d\mu}^2\leq\; \mu(\{f\geq t_0\})\int g^{2q}\,d\mu\leq\frac{1}{2}\int g^{2q}d\mu\;,
$$
$$
\var{\mu}{g^q}\geq \frac{1}{2}\int g^{2q}\,d\mu\geq \frac{1}{2}\,\mu(\{f\geq 2t\})=\frac 12\,a(2t)\;.
$$
On the other hand, by the $L^q$-Poincar\'e inequality, 
$$
\PAR{\frac 12\,a(2t)}^{1/q}\!\leq\;\PAR{\frac 12\,\var{\mu}{g^q}}^{1/q}\leq\; \CP\!\int|\nabla g|^2\,d\mu\;\leq\; \CP\,\frac{\mu(\{t\leq \rho<2t\})}{t^2}
$$
using the fact that, a.e., $\ABS{\nabla f}^2\leq \NRM{f}_{\rm{Lip}}\leq 1$. With $\kappa:=2^{1/q}\,\CP$, this proves that 
$$
t^2\,\PAR{a(2t)}^{1/q}+\kappa\,a(2t)\leq \kappa\,a(t)\quad\forall\; t>t_0\;,
$$
and as a consequence, 
$$
t^2\,\PAR{a(2t)}^{1/q}\leq \kappa\,a(t)\quad\forall\; t>t_0\;.
$$
With $c:=\kappa^q\,t_0^{-2q}$ and $a_n:=a(2^n\,t_0)$, this means
$$
a_{n+1}\leq c\,2^{-2nq}\,a_n^q\quad\forall\;n\in\dN\;.
$$
If $b_n:=2^{2n\,\alpha_n}$ with $(n+1)\,\alpha_{n+1}=n\,q\,\alpha_n-2\,n\,q$, then $b_{n+1}\leq c\,b_n^q$ and $\limsup_{n\to\infty}b_n$ is therefore bounded by the unique fixed point, $\bar b$, of $b\mapsto c\,b^q$. The sequence $(\alpha_n)_{n\in\dN}$ converges to $q/(1-q)$. Hence
$$
a(t)\leq a_n\leq O\PAR{2^{-2n\,\frac q{1-q}}}\quad\mbox{as}\quad t\to\infty\;,
$$
where $n$ is the integer part of $\log(t/t_0)/\log 2$. This concludes the proof. 
\end{eproof}

\section{Application to the weighted porous media equation}
\label{sec-app}

Let $d$ be a positive integer and $\psi\in\Cn{2}(\dR^d)$ a function such that $\int e^{-\psi}dx<+\infty$. We define the probability measure
$$
d\mu_\psi:=\frac{e^{-\psi}\,dx}{Z_\psi}
$$
and the operator $\GI$ on $\Cn{2}(\dR^d)$ by
$$
\forall\;f\in \Cn{2}(\dR^d)\,,\quad\GI f:=\Delta f-\nabla \psi\cdot\nabla f\;.
$$
Such a generator $\GI$ is symmetric in $L^2_{\mu_\psi}(\dR^d)$,
$$
\forall\;f,\,g\in\Cn{1}(\dR^d)\,,\quad\int f\,\GI g\,d\mu_\psi=-\int\nabla f\!\cdot\!\nabla g\,d\mu_\psi\;.
$$
We consider for $m>1$ the nonlinear partial differential equation
\begin{equation*}
\label{eq-pde}
\frac{\partial u}{\partial t}=\GI\,u^m\eqno\WPME
\end{equation*}
for $t\geq 0$, $x\in\dR^d$, corresponding to a non-negative initial condition $u(0,x)=u_0(x)$ for any $x\in\dR^d$. Such an equation will be called the {\sl weighted porous media equation.\/}

\subsection{$L^1$-contraction, existence and uniqueness}
\label{sec-exi}

The existence proof is based on the method developed by V\'azquez in \cite{vazquez}. The main difference between the standard porous media and the weighted porous media equations is that a natural space to study weak solutions of \WPME is a weighted space, for instance $L^2_{\mu_\psi}(\dR^d)$, which contains all constant functions. We shall first consider the case of a bounded domain and then extend solutions to the whole space.\medskip

Consider first a bounded domain $\Omega\subset\dR^d$ with smooth boundary. Denote by $Q=\Omega\times [0,+\infty)$, $\Sigma=\partial\Omega\times [0,+\infty)$. Let $u_0$ be a positive function in $\Omega$ which satisfies $n\cdot\nabla u_0=0$ on $\partial \Omega$, where $n=n(x)$ denotes the outgoing normal unit vector at $x\in\partial\Omega$. We shall say that $u$ is classical solution of \WPME in $\Omega$ if $u$ is a $\Cn{2}$ function on $Q$ such that
\begin{equation}\label{Eqn:WPME-bded}
\left\{
\begin{array}{rl}
\disp u_t=\GI\,u^m&\rm{ in }\quad Q\\
\disp u(\cdot,0)=u_0&\rm{ in } \quad\Omega\\
\disp n\cdot\nabla u=0&\rm{ on } \quad\Sigma
\end{array}
\right.
\end{equation}
\begin{elem}[$L^1$-contraction principle]
\label{lem-contrac}
Let $\Omega\subset\dR^d$ be a bounded domain with smooth boundary. If $u$, $\hat{u}$ are two classical solutions of \eqref{Eqn:WPME-bded} with smooth positive initial data $u_0$ and~$\hat{u}_0$, then for all $t>\tau>0$ one gets
$$
\int_\Omega (u(t,x)-\hat{u}(t,x))_+\,d\mu_\psi(x)\leq\int_\Omega (u(\tau,x)-\hat{u}(\tau,x))_+\,d\mu_\psi (x)\leq\int_\Omega (u_0(x)-\hat{u}_0(x))_+\,d\mu_\psi(x)\;.
$$
\end{elem}
\begin{eproof}
Let $\chi\in\Cn{1}(\dR)$ be such that $0\leq\chi\leq 1$, $\chi(s)=0$ for $s\leq0$, $\chi'(s)>0$ for $s>0$. We obtain
\begin{eqnarray*}
\int_\Omega (u-\hat{u})_t\,\chi(u^m-\hat{u}^m)\,d\mu_\psi&=&\int_\Omega \GI(u^m-\hat{u}^m)\,\chi(u^m-\hat{u}^m)\,d\mu_\psi\\
&=&-\int_\Omega |\nabla(u^m-\hat{u}^m)|^2\,\chi'(u^m-\hat{u}^m)\,d\mu_\psi\;,
\end{eqnarray*}
using that $n\cdot\nabla(u^m-\hat{u}^m)=0$ on $\Sigma$. Therefore, by taking $\chi$ as a smooth approximation of the function ${\rm sgn}^+_0$ which is identically equal to $1$ on $(0,+\infty)$ and to $0$ on $(-\infty,0]$, and observing that
$$
\frac\partial{\partial t}(u-\hat{u})_+=\mbox{\rm sgn}^+_0(u-\hat{u})\,\frac\partial{\partial t}(u-\hat{u})\;,
$$
we get that
$$
\frac{d}{dt}\int_\Omega (u-\hat{u})_+\,d\mu_\psi\leq0\;.
$$
\end{eproof}

Lemma~\ref{lem-contrac} results in a Maximum Principle for (WPME).
\begin{ecor}[Maximum Principle and uniqueness]
\label{cor-max}
Let $\Omega\subset\dR^d$ be a bounded domain with smooth boundary and consider two classical solutions $u$, $\hat{u}$ of \eqref{Eqn:WPME-bded} with smooth positive initial data $u_0$ and $\hat{u}_0$. If $u_0\leq\hat{u}_0$ in $\Omega$, then $u\leq\hat{u}$ in $Q$. As a consequence, the classical solution of \WPME is unique.
\end{ecor}
We may now apply the existence theory for non degenerate parabolic equations as in \cite{lsu}.
\begin{eprop}
\label{prop-ex1}
Let $\Omega\subset\dR^d$ be a bounded domain with smooth boundary. For any positive function $u_0\in\Cn{\infty}$, there exits a classical solution $u$ of \eqref{Eqn:WPME-bded} with initial datum $u_0$. Moreover, for all $\tau\geq0$,
\begin{equation*}
\label{eq-reg}
\int_\Omega\int_0^\tau|\nabla u(t,x)^m|^2\,dt\,d\mu_\psi(x)-\int_\Omega u^{m+1}(\tau,x)\,d\mu_\psi(x)=-\int_\Omega u_0(x)^{m+1}\,d\mu_\psi(x)\;.
\end{equation*}
\end{eprop}
\begin{eproof}
Since the initial datum $u_0$ is positive on $\bar{\Omega}$, which is bounded, it follows that $\min(u_0)>0$. Standard quasilinear theory (see chapter 6 of \cite{lsu}) applies, thus providing us with a classical solution of (WPME). For all $\tau>0$, we have
\begin{eqnarray*}
\int_\Omega u^{m+1}(\tau,x)\,d\mu_\psi(x)-\int_\Omega u_0(x)^{m+1}\,d\mu_\psi(x)&=&\int_\Omega\left(\int_0^\tau\frac{\partial}{\partial s}\,u^{m+1}(s,x)\,dt\right)d\mu_\psi(x)\\
&=&-(m+1)\int_\Omega\int_0^\tau|\nabla u^m|^2\,d\mu_\psi\,dt
\end{eqnarray*}
after an integration by parts.
\end{eproof}
\medskip

The results obtained for a bounded domain can be extended to solutions in the whole euclidean space. Various results can be stated which are out of the scope of our paper, so let us make some simplifying assumptions.
\begin{eprop}
\label{cor-ex}
Let $u_0$ be a $\Cn{1}\cap L^{m+1}_{\mu_\psi}(\dR^d)$ positive initial condition. Then there exists a unique classical solution of the \WPME with initial datum $u_0$.
\end{eprop}
\begin{eproof} We can approximate the solution using the following scheme:
\begin{enumerate}[{\rm (1)}]
\item Consider an initial datum which is uniformly bounded away from $0$, for instance $u_0^n=u_0+1/n$. 
\item Consider a regularized drift term $\psi_n$ such that $\psi_n\equiv\psi$ in $B(0,n)$ and $\psi_n(x)\equiv c_n\,|x|^2/2$ in $B(0,n)^c$, with $c_n>0$. 
\item Modify the operator $\GI$ by considering $\GI_nf^m:=\nabla(a_n(f)\,\nabla f)-\nabla\psi_n\cdot\nabla f^m$ where $a_0$ is a smooth positive function on $[0,m\,n^{-(m-1)}]$ and such that $a_n(s)=m\,s^{m-1}$ for any $s\geq 1/n$. The standard theory of parabolic equations applies and provides us with a solution $u_n$ of the regularized equation, ${\partial u_n}/{\partial t}=\GI_n\,u_n^m$. 
\item Prove an $L^1$-contraction principle as in Lemma~\ref{lem-contrac}, from which we deduce a Maximum Principle and the uniqueness of the approximating solution $u_n$. We observe that $u_n$ is a solution of (WPME) with initial datum $u_0^n$, except that $\psi$ has to be replaced by $\psi_n$. 
\item Barrier functions based on the solution of the heat equation can be provided, thus showing the conservation of the $L^1$-norm (with respect to the measure $d\mu_\psi$) and uniform estimates with respect to $n\in\dN$.
\item Take a pointwise monotone limit as $n\to\infty$ and obtain a weak solution $u(t,x)$ of (WPME) with initial datum $u_0$. Classical regularity properties (see for instance \cite[Chap. 6]{lsu}) prove that the weak solution is a classical solution on $\dR^d$.
\end{enumerate}
\end{eproof}

\subsection{Asymptotic behavior of the solutions}
\label{sec-asy}

\begin{ethm}
Let $m\geq1$.
\begin{enumerate}
\item [{\rm (i)}] If $(\mu_\psi,\mu_\psi)$ satisfies a $L^q$-Poincar\'e inequality, $q=2/(m+1)$, for some constant $\CP >0$, then for any initial condition $u_0\in L^2(\mu_\psi)$, we have
$$
\forall\;t\geq0\;,\quad\var{\mu_\psi}{u(\cdot,t)}\leq\PAR{\big[\var{\mu_\psi}{u_0}\big]^{-(m-1)/2}+ \frac{4\,m\,(m-1)}{(m+1)^2}\,\CP\,t}^{-2/(m-1)}\;.
$$
Reciprocally, if the above inequality is satisfied for any $u_0$, then $(\mu_\psi,\mu_\psi)$ satisfies a $L^q$-Poincar\'e inequality with constant $\CP$.
\item [{\rm (ii)}] If $(\mu_\psi,\mu_\psi)$ satisfies a $L^q$-logarithmic Sobolev inequality, $q=1/m$, for some constant $\CLSI >0$, then for any non-negative initial condition $u_0$ such that $\ent{\mu_\psi}{u_0}<\infty$, we have
$$
\forall\;t\geq0\;,\quad\ent{\mu_\psi}{u(\cdot,t)}\leq\PAR{\big[\ent{\mu_\psi}{u_0}\big]^{1-m}+\frac{4\,(m-1)}{m}\,\CLSI\,t}^{-1/(m-1)}\;.
$$
Reciprocally, if the above inequality is satisfied for any $u_0$, then $(\mu_\psi,\mu_\psi)$ satisfies a $L^q$-logarithmic Sobolev inequality with constant $\CLSI$.
\end{enumerate}
\end{ethm}
\begin{eproof}
Let us briefly sketch the first result.
\begin{eqnarray*}
\frac d{dt}\,\var{\mu_\psi}{u}=2\int u_t\,u\,d\mu_\psi
=2\int u \,\GI u^m \,d\mu_\psi
&=&-2\int\nabla u\cdot\nabla (u^m)\,d\mu_\psi\\
&=&-\frac{8m}{(m+1)^2}\int|\nabla u^{\frac{m+1}2}|^2\,d\mu_\psi\;.
\end{eqnarray*}
One can now apply the $L^q$-Poincar\'e inequality with $u=f^{2/(m+1)}$, $q=2/(m+1)$, to get
$$
\frac d{dt}\,\var{\mu_\psi}{u}\le -\frac{8m}{(m+1)^2}\,\CP\left[\var{\mu_\psi}{u}\right]^{\frac{m+1}2}\;.
$$
A simple integration of this differential inequality gives the result. Reciprocally, a derivation at $t=0$ gives the $L^q$-Poincar\'e inequality with constant $\CP$. The proof in the second case is similar.
\end{eproof}
\begin{exe}
Consider on $\dR^d$ the probability measure $\mu$ given by
$$
d\mu(x_1,\cdots,x_n):=Z^{-1}\,e^{W(x_1,\cdots,x_n)}\,\prod_{i=1}^n\frac{dx_i}{\PAR{1+\ABS{x_i}}^{1+\al}}\;,
$$
where $Z$ is a normalization constant and $W$ a bounded function on $\dR^d$. By Proposition~\ref{Prop:Pertub-Tensor-Concentr}, the measure $\mu$ satisfies a $L^q$-Poincar\'e inequality with $q\in[1/2,1)$ if $\al>2q/(1-q)$. Then the variance of the solution to the associated \WPME converges to $0$ as $t\to+\infty$ if $m>(\al+4)/\al$.
\end{exe}

\section{Appendix. A Variant of two results of \cite{ca-ba-ro2}}
\label{Sec:Appendix}

We present variants of Theorem 1 and 2 in \cite{ca-ba-ro2}, in which we remove any assumption on the absolute continuity of the measure $\mu$ with respect to the volume measure.

\medskip
We recall that $(\mu,\nu)$ satisfies a {\sl weak Poincar\'e inequality\/} with associated function $\CWP$ if
$$
\var{\mu}f\leq\CWP\PAR{s}\int\ABS{\nabla f}^2 \,d\nu+s\,\big[\osc{f}\big]^2\quad\forall\;s\in (0,1/4)\,,\quad\forall\;f\in{\mathcal C}^1(M)\;.
$$
See Definition~\ref{Defn:WeakPoincare} for details.
\begin{ethm}{\rm \cite{ca-ba-ro2}}
Let $\mu$ and $\nu$ be respectively a probability measure and a positive measure on~$M$. Assume that that $(\mu,\nu)$ satisfies a weak Poincar\'e inequality for some non-negative non increasing function $\CWP(s)$. Let $\gamma(s):=4\,\CWP(s/4)$. Then for every measurable subsets $A$, $B$ of $M$ such that $A\subset B$ and $\mu(B)\leq 1/2$,
$$
\Capa_\nu(A,B)\ge \frac{\mu(A)}{\gamma(\mu(A))}\;.
$$
\end{ethm}
\begin{eproof} The proof of \cite{ca-ba-ro2} can be extended to the case of two measures $\mu$ and $\nu$ without changes. Let us sketch it for completeness. Let $f$ be such that $\1_A\leq f\leq \1_B$ and observe that $\osc f\leq 1$. By the Cauchy-Schwarz inequality, $\PAR{\int f\,d\mu}^2\leq\mu(B)\int f^2\,d\mu\leq \frac 12\int f^2\,d\mu$. Hence
$$
\CWP\PAR{s}\int\ABS{\nabla f}^2 \,d\nu+s\geq\var{\mu}f\ge \frac 12\int f^2\,d\mu\ge\frac{\mu(A)}2\;,
$$
which completes the proof after noticing that $\frac a{\gamma(a)}=\frac a{4\,\CWP(a/4)}\leq\sup_{s\in(0,1/4)}\frac{a/2-s}{\CWP(s)}$ with $a/2=\mu(A)/2\leq 1/4$.
\end{eproof}

\medskip
In the next result, we explicitly remove the assumption of absolute continuity with respect to the volume measure. Let $\kkappa:=(11+5\sqrt 5)/2\approx 11.0902$.
\begin{elem}
\label{prel-crit-capa}
Let $\mu$ and $\nu$ be respectively a probability measure and a positive measure on~$M$. For some $\theta\in (0,1)$, consider a positive non increasing function $\gamma$ on $(0,\theta)$ and assume that for every measurable subsets $A$, $B$ of $M$ such that $A\subset B$ and $\mu(B)\leq\theta$,
$$
\Capa_\nu(A,B)\ge \frac{\mu(A)}{\gamma(\mu(A))}\;.
$$
Then for every function $f\in\Cn{1}(M)$ such that $\mu(\Omega_+)\leq\theta$, $\Omega_+:=\{f>0\}$, and every $s\in(0,1)$ one has
$$
\int f_+^2 \le \kkappa\,\gamma(s) \int_{\Omega_+}|\nabla f|^2\,d\nu+ s\, \Big[\mathop{\mbox{\rm supess}_\mu}_{}f\Big]^2\;.
$$
\end{elem}
\begin{eproof}
Fix $s\in (0,1)$. Let $c=c(s):=\inf\{t\ge 0\,:\,\mu(\{f_+>t\})\le s \}$. If $c=0$ then $\mu(\Omega_+)\le s$ and $\int_{\Omega_+} f_+^2\,d\mu \le s\,\mbox{\rm supess}_\mu f_+^2$. We also know that $c\leq\mbox{\rm supess}_\mu f_+$ and $s_0:=\mu(\{f_+>c\})\le s$. For a given $\rho\in (0,1)$, let $\Omega_k:=\{f_+>c\,\rho^k\}$ for any $k\in\mathbb N$ and define $s_k:=\mu(\Omega_k)$. We observe that
$$
s_0\leq s\leq s_1\leq s_2\leq\ldots\leq \mu(\Omega_+)\leq\theta\;.
$$
Using the decomposition $\Omega_+=\Omega_0+\bigcup_{k\in\mathbb N}(\Omega_{k+1}\setminus\Omega_k)$, we get
$$
\int_{\Omega_+}f_+^2\,d\mu=\int_{\Omega_0}f_+^2\,d\mu+\sum_{k>0}\int_{\Omega_{k+1}\setminus\Omega_k}f_+^2\,d\mu\leq s_0 \,\SBRA{\mbox{\rm supess}_\mu f_+}^2 \;+\; \sum_{k\in\mathbb N} c^2\,\rho^{2k}\,(s_{k+1}-s_k)\;.
$$
Since $s_0\leq s$, we can actually write
$$
\int_{\Omega_+}f_+^2\,d\mu=\int_{\Omega_0}f_+^2\,d\mu+\sum_{k>0}\int_{\Omega_{k+1}\setminus\Omega_k}f_+^2\,d\mu\leq s\,\SBRA{\mbox{\rm supess}_\mu f_+}^2 \;+\; \sum_{k\in\mathbb N} c^2\,\rho^{2k}\,(s_{k+1}-s_k)\;.
$$
On the other hand, by noticing that
$$
\sum_{k\in\mathbb N}\rho^{2k}\,(s_{k+1}-s_k)=\left(\frac 1{\rho^2}-1\right)\sum_{k\in\mathbb N}\rho^{2k}\,s_k-\frac 1{\rho^2}\,s_0=\frac{1-\rho^2}{\rho^2}\sum_{k\in\mathbb N}\rho^{2k}\,(s_k-s_0)\;,
$$
we get
$$
\int_{\Omega_+}f_+^2\,d\mu \le s\,\SBRA{\mbox{\rm supess}_\mu f_+}^2\;+\;c^2\,\frac{1-\rho^2}{\rho^2}\sum_{k\in\mathbb N}\rho^{2k}\,(s_k-s_0)\;.
$$
By our assumptions, for any $k\in\mathbb N\setminus\{0\}$, the function $\theta\mapsto\frac{s_k-\,s_0+\,\theta\,s_0}{\gamma(s+\theta\,(s_k-s))}$ is monotone increasing on $(0,1)$. Hence
$$
\frac{s_k-s_0}{\gamma(s)}\leq\frac{s_k}{\gamma(s_k)}\leq\Capa_\nu(\Omega_k,\Omega_{k+1})\;,
$$
$$
s_k-s_0\le \gamma(s)\,\Capa_\nu(\Omega_k,\Omega_{k+1})\;.
$$
Define
$$
F_k:=\min\left\{1,\frac 1{1-\rho}\left( \frac{f_+}{c\,\rho^k}-\rho\right)_{\!\!+}\,\right\}\,.
$$
By definition of $\Capa_\nu(\Omega_k,\Omega_{k+1})$, we have
$$
\Capa_\nu(\Omega_k,\Omega_{k+1})\leq\int_{\Omega_{k+1}\setminus\Omega_k}|\nabla F_k|^2\,d\mu=\frac 1{c^2\,\rho^{2k}\,(1-\rho)^2}\int_{\Omega_{k+1}\setminus\Omega_k}|\nabla f_+|^2\,d\mu\;.
$$
Collecting the estimates, we get
$$
\int_{\Omega_+}f_+^2\,d\mu \le s\,\SBRA{\mbox{\rm supess}_\mu f_+}^2\;+\;\frac 1{\rho^2}\,\frac{1+\rho}{1-\rho}\int_{\Omega_+}|\nabla f_+|^2\,d\mu\;.
$$
The result follows after optimizing on $\rho$, that is by taking $\rho=\PAR{\sqrt 5-1}/2$.
\end{eproof}

\medskip
Theorem 2 in \cite{ca-ba-ro2} can be generalized in case of two measures $\mu$ and $\nu$ which are not absolutely continuous with respect to the volume measure as follows. The main idea is to apply Lemma~\ref{prel-crit-capa} with $\theta=1/2$ and use the median of $f$ to define subsets of $M$ with measure at most equal to $\theta$.
\begin{ethm}
\label{crit-capa}
Let $\mu$ and $\nu$ be respectively a probability measure and a positive measure on~$M$. Consider a positive non increasing function $\gamma$ on $(0,1/2)$ and assume that for every measurable subsets $A$, $B$ of $M$ such that $A\subset B$ and $\mu(B)\leq 1/2$,
$$
\Capa_\nu(A,B)\ge \frac{\mu(A)}{\gamma(\mu(A))}\;.
$$
Then with $\kkappa:=(11+5\sqrt 5)/\,2$, for every function $f\in\Cn{1}(M)$ and every $s\in(0,1/4)$ one has
$$
\var\mu f \le \kkappa\,\gamma(s) \int |\nabla f|^2\,d\nu+ s\,\big[\osc{f}\big]\;.
$$
\end{ethm}
\begin{eproof}
Fix $s\in (0,1/4)$ and let $m$ be a median of $f$ with respect to $\mu$. Denote by $\Omega_+$ and~$\Omega_-$ the sets $\{f>m\}$ and $\{f<m\}$. By definition of $m$, $\mu(\Omega_\pm)\leq 1/2$. By definition of the variance, we also know that
$$
\var\mu f\leq\int(f-m)^2\,d\mu=\int_{\Omega_+} (f-m)^2\,d\mu+\int_{\Omega_-} (f-m)^2\,d\mu\;.
$$
We can apply Lemma~\ref{prel-crit-capa} to $F_+$ and $F_-$ with $\theta=1/2$ and get the result using the fact that, if $a={\rm infess}_\mu f$ and $b={\rm supess}_\mu f$, then for any $m\in[a,b]$, $(b-m)^2+(m-a)^2\leq(b-a)^2=\osc f^2$.
\end{eproof}

\medskip
Concerning the weak logarithmic Sobolev inequalities, the absolute continuity of $\mu$ with respect to the volume measure can easily be removed in Theorem 2.1 of \cite{ca-ge-gu}, without change in the proof.
\begin{elem}\label{Lem:CGG05-2.2}
Let $\mu$ and $\nu$ be respectively a probability measure and a positive measure on~$M$. If $(\mu,\nu)$ satisfies a weak logarithmic Sobolev inequality with associated function $\CWL$, then for every $A\subset B\subset M$ such that $\mu(B)\leq 1/2$,
$$
\Capa_\nu(A,B)\ge\frac{\mu(A)\log\left(1+\frac{e^2}{\mu(A)}\right)}{\gamma\SBRA{\mu(A)\log\left(1+\frac{e^2}{\mu(A)}\right)}}
$$
with $\gamma(s):= 2\,\CWL(s/2)$.
\end{elem}
\begin{eproof}
As in the proof of Proposition~\ref{cor-1ls}, it holds that
$$
\Capa_\nu(A,B)\ge\frac{\mu(A)\log\left(1+{e^2}/{\mu(A)}\right)-s}{\CWL(s)}\;.
$$
\end{eproof}

\medskip
Reciprocally, if the capacity measure criterion of Lemma~\ref{Lem:CGG05-2.2} is satisfied, it is not clear that a weak logarithmic Sobolev inequality holds unless we assume the absolute continuity of $\mu$ with respect to the volume measure. See Theorem 2.2 of \cite{ca-ge-gu} in such a case.

\bigskip\begin{spacing}{0.8}{\small {\bf Acknowledgments.} {\sl The authors thank Philippe Lauren\c cot for helpful discussions on existence results for weighted porous media equations. They have been supported by the IFO project of the French Research Agency (ANR).}}\end{spacing}
\newcommand{\etalchar}[1]{$^{#1}$}
\def\cprime{$'$}

\bigskip J. Dolbeault and I. Gentil: CEREMADE, Univ. Paris Dauphine, Place de Lattre de Tassigny, F-75775 Paris Cedex 16, France. \email{$\{$dolbeaul, gentil$\}$@ceremade.dauphine.fr}

\smallskip A. Guillin: Ecole Centrale de Marseille  and Laboratoire d'Analyse, Topologie, Probabilit\'es, Univ. de Provence, 39, rue F. Joliot-Curie, F-13453 Marseille Cedex 13, France. \email{guillin@ceremade.dauphine.fr}

\smallskip F.Y. Wang: School of Mathematical Sciences, Beijing Normal Univ., Beijing 100875, China. \email{wangfy@bnu.edu.cn}

\end{document}